\documentclass[12pt,a4paper]{article}%% ,draft
\usepackage{amssymb,ngerman,graphicx}
\newcounter{zahl}
\newenvironment{publikationen}[1]{%
\begin{list}%
{\rm\small\arabic{zahl}.}{\usecounter{zahl}%
\setlength\leftmargin\leftmarginiv%
\settowidth\labelwidth{#1.}%
\setlength\labelsep\leftmargin%
\addtolength\labelsep{-\labelwidth}%
\small\itemsep0.3em plus 0.1em minus 0.1em
\parsep0em plus 0.1em \rm}}{\end{list}}

\begin{document}
\title{Erinnerungen an Heinrich Brauner (1928--1990)%
\thanks{Ausarbeitung eines Vortrags beim 33.\ S\"{u}ddeutschen Differentialgeometrie-Kolloquium
am 23.\ Mai 2008 an der Technischen Universit\"{a}t Wien.}}
\author{Hans Havlicek}
\date{}

\maketitle

%%%%%%%%%%%%%%%%%%%%%%%%%%%%%%%
\begin{abstract}
Zur achtzigsten Wiederkehr des Geburtstages von Heinrich Brauner sollen
Eigenschaften, Wesensmerkmale und Leistungen dieses \"{o}sterreichischen Geometers
aufgezeigt werden. Dabei m\"{o}chte ich zumindest ein wenig von dem vermitteln, was
aus meiner Sicht das Besondere dieses au{\ss}ergew\"{o}hnlichen Wissenschaftlers und
Menschen ausmachte.


\end{abstract}
%%%%%%%%%%%%%%%%%%%%%%%%%%%%%%%%%

\section{Einleitung}\label{sec:pre}

In meinem Vortrag m\"{o}chte ich einige Worte der pers\"{o}nlichen Erinnerung an meinen
Lehrer,
\begin{center}
    Herrn O.Univ.Prof.\ Mag.rer.nat.\ Dr.phil.\ Dr.techn.\ Heinrich Brauner,
\end{center}
sprechen, dessen Geburtstag sich im Jahr 2008 zum achtzigsten Male j\"{a}hrt.
Zun\"{a}chst sei sein Lebensweg ganz kurz skizziert, wobei ich mich auf die Angaben
in \cite{1} st\"{u}tze.
\par
\begin{figure}[!t]
\centering
    \includegraphics[width=0.9\textwidth]{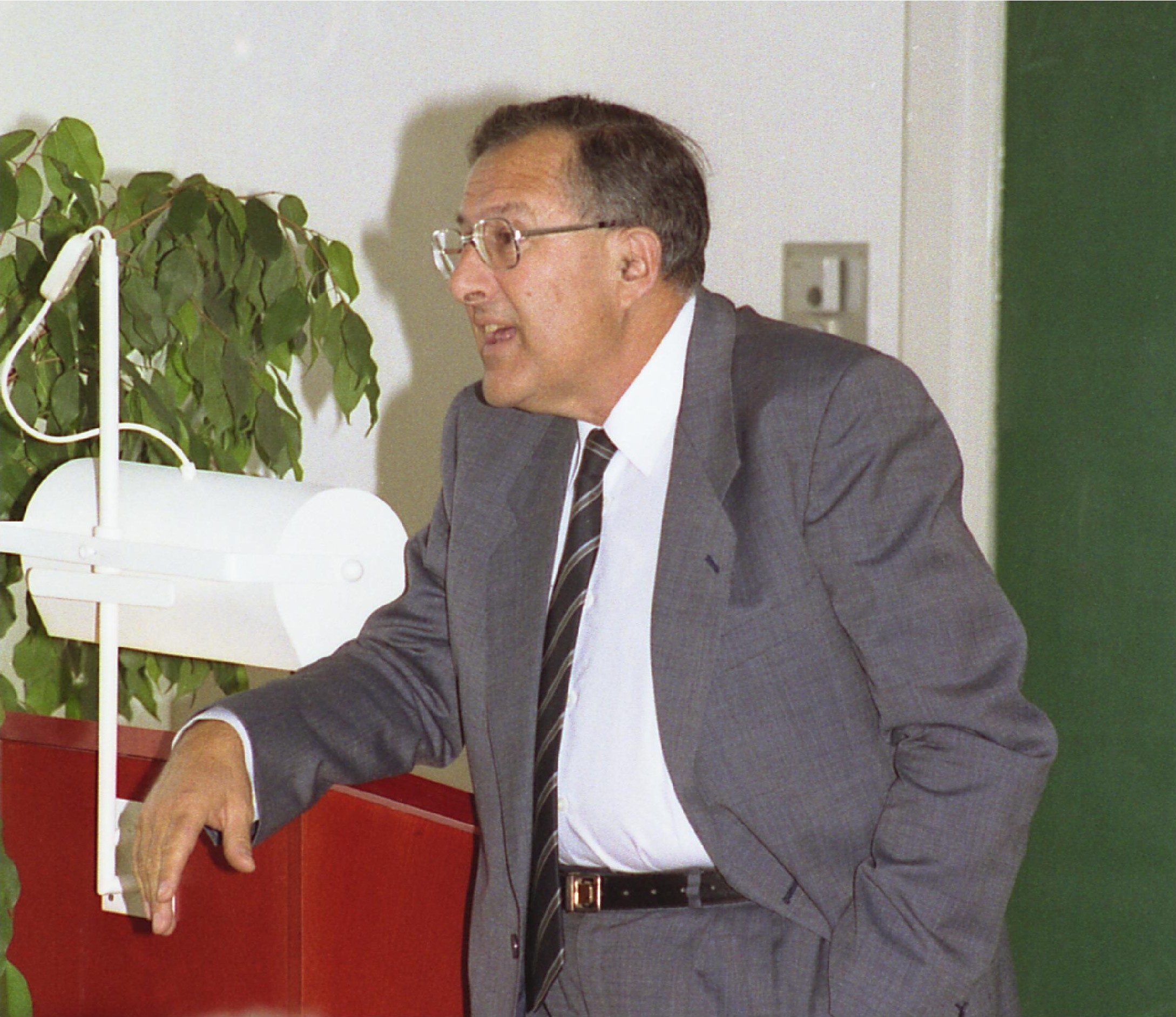}
    \caption{Festkolloquium 1988}\label{bild1} %% height=9cm
\end{figure}
Heinrich Brauner wurde am 21.~November 1928 in Wien geboren, wo er auch das
Realgymnasium besuchte. Von 1946 bis 1952 studierte er an der Universit\"{a}t Wien
und der Technischen Hochschule Wien. Brauner legte die Lehramtspr\"{u}fungen f\"{u}r
die F\"{a}cher \emph{Mathematik}, \emph{Physik} und \emph{Darstellende Geometrie}
und die erste Staatspr\"{u}fung aus \emph {Technischer Physik} ab. Er verfasste
zwei Dissertationen: \emph{\"{U}ber $n+1$ fache Orthogonalsysteme von Riemannschen
Hyperfl\"{a}chen der Klasse $1$ im euklidischen Raum $R^{n+1}$} bei Johann Radon
sowie \emph{Kongruente Verlagerung kollinearer R\"{a}ume in axiale Lage} bei Walter
Wunderlich. Brauner wurde an der Universit\"{a}t Wien zum Doktor der Philosophie
und an der Technischen Hochschule Wien zum Doktor der Technischen
Wissenschaften promoviert.
\par
\begin{figure}[!t]
\centering
    \includegraphics[width=0.9\textwidth]{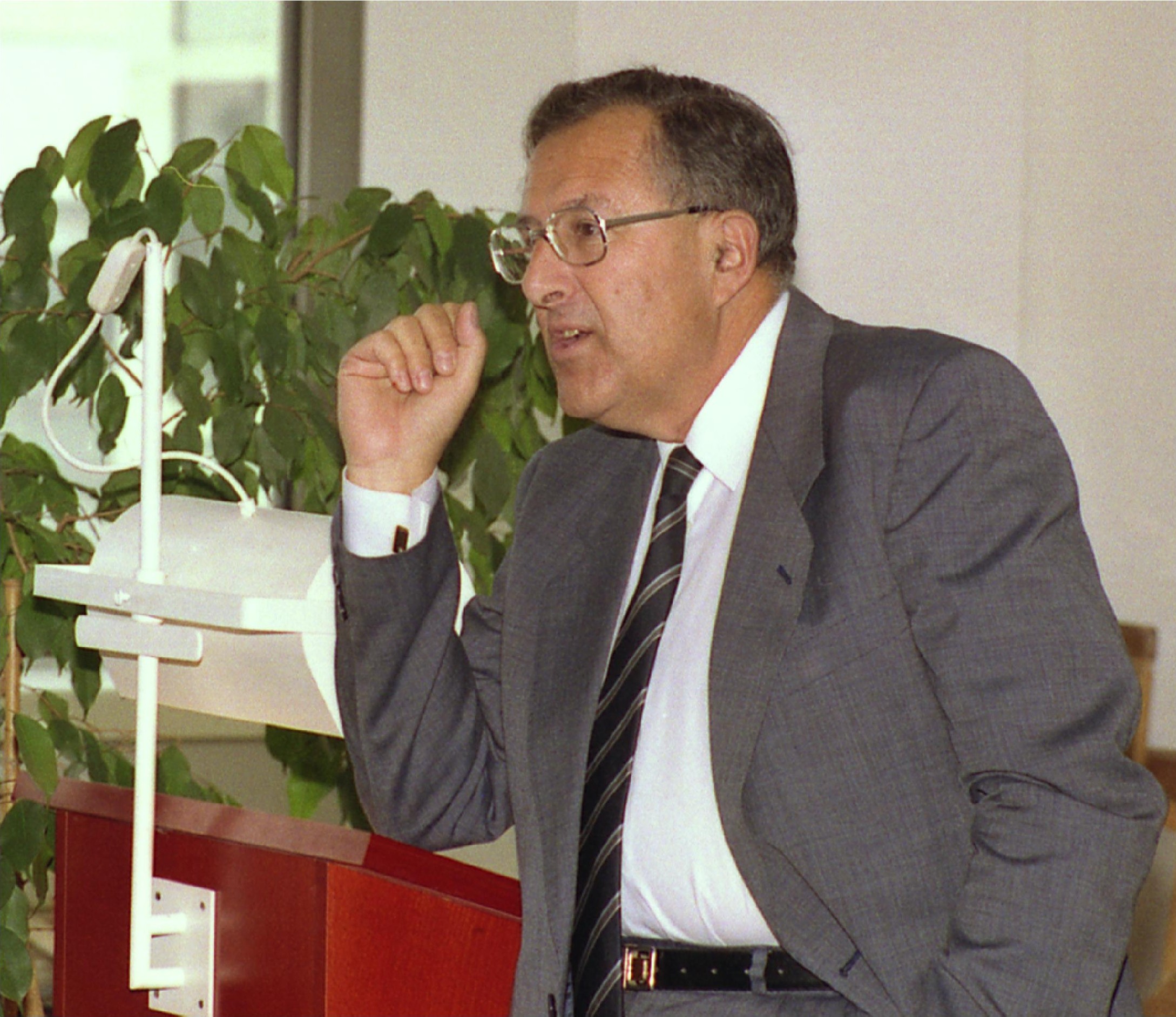}
    \caption{Festkolloquium 1988}\label{bild2} %% height=9cm
\end{figure}
Ab 1950 war Brauner im Schuldienst t\"{a}tig und daneben ab 1951 teilbesch\"{a}ftigte
wissenschaftliche Hilfskraft am 1.\ Institut f\"{u}r Geometrie der Technischen
Hochschule Wien. Erst 1954 konnte er ebendort eine Stelle als vollbesch\"{a}ftigter
Hochschulassistent antreten.
\par
Schon 1956 habilitierte sich Brauner an der Technischen Hochschule Wien f\"{u}r das
Fach \emph{Geometrie, insbesondere Darstellende Geometrie\/} und im Jahr
darauf, in einem davon unabh\"{a}ngigen Verfahren, an der Universit\"{a}t Wien f\"{u}r das
Fach \emph{Mathematik}.
\par
Im Jahre 1960 nahm er einen Ruf auf ein Ordinariat an der Technischen
Hochschule Stuttgart an. Ab 1969 war Brauner Ordentlicher Universit\"{a}tsprofessor
f\"{u}r Geometrie an der Technischen Hochschule (Technischen Universit\"{a}t) Wien.
\par %% Akademie seit 16.5.1972
Brauner wurde im Jahr 1970 zum Honorarprofessor der Universit\"{a}t Wien ernannt.
Ferner war er ab 1972 korrespondierendes Mitglied der \"{O}sterreichischen Akademie
der Wissenschaften und in weiterer Folge Tr\"{a}ger des Ehrenkreuzes f\"{u}r
Wissenschaft und Kunst I.~Klasse.

\par
Die Abbildungen \ref{bild1} und \ref{bild2} zeigen Brauner beim Festkolloquium,
das aus Anlass seines 60.\ Geburtstages am 21.\ Oktober 1988 am Institut f\"{u}r
Geometrie der Technischen Universit\"{a}t Wien stattfand. Er litt zu diesem
Zeitpunkt bereits an Osteoporose. Brauner k\"{a}mpfte gegen diese sehr schmerzhafte
Krankheit mit unendlicher Geduld an und nahm seine Aufgaben am Institut bis
wenige Wochen vor seinem Ableben wahr. Heinrich Brauner erlag seinem schweren
Leiden am 1.~Juni 1990.

\section{Der Lehrer Heinrich Brauner: Es begann mit einem Punktsack}

Meine erste Begegnung mit Heinrich Brauner war im Wintersemester 1972/73 in
seiner Vorlesung \emph{Projektive Geometrie I} f\"{u}r die erstj\"{a}hrigen
Lehramtskandidaten. Ich hatte keine Ahnung, was mich erwarten w\"{u}rde -- weder
inhaltlich noch den Vortragenden betreffend.
\begin{figure}[!h]
\centering
    \includegraphics[width=10.5cm]{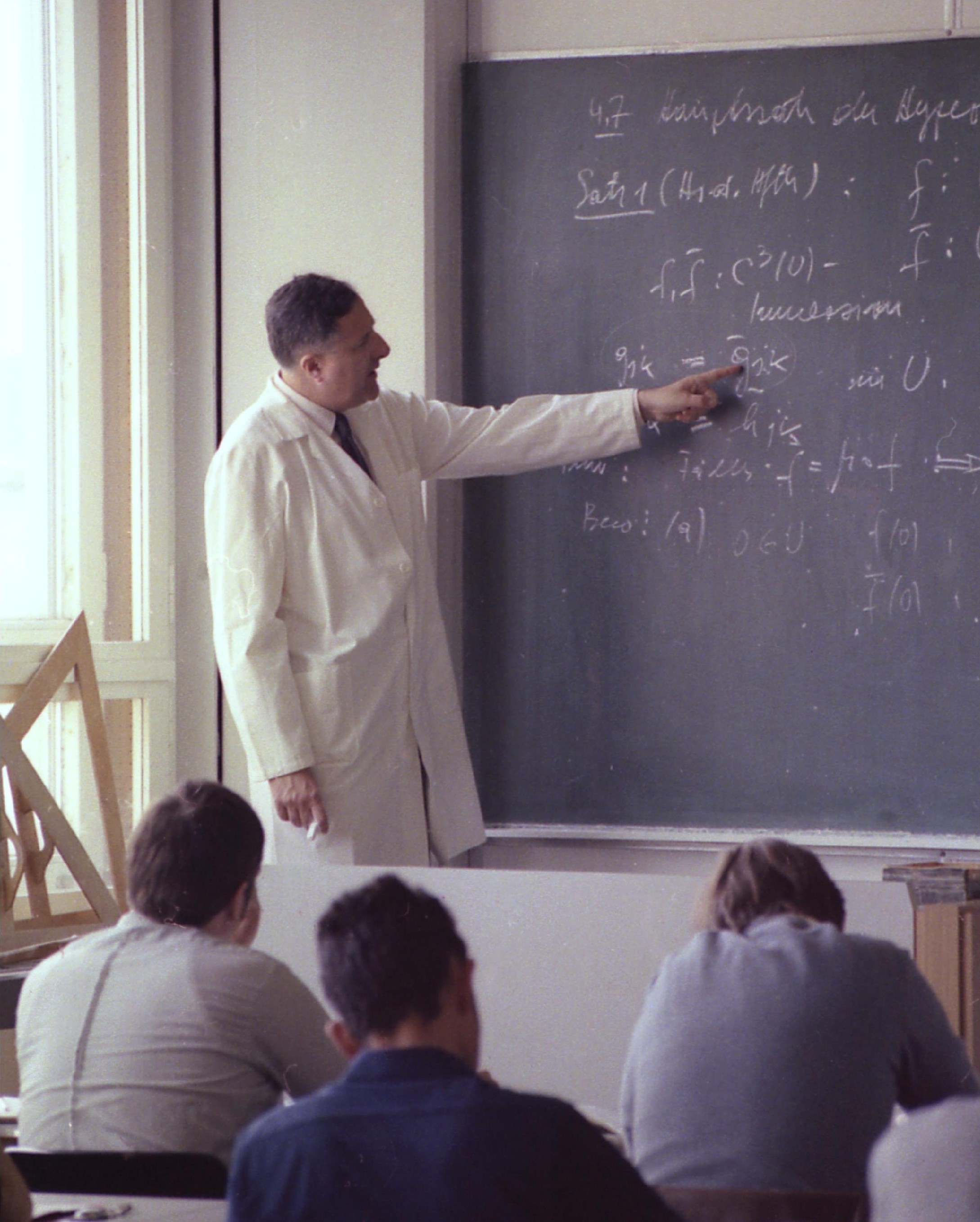}
    \caption{Vorlesung \"{u}ber Differentialgeometrie 1982}\label{bild3}
\end{figure}
Brauner begann die erste Vorlesung und brachte sogleich mein in der Schule
erworbenes Bild der Geometrie kr\"{a}ftig ins Wanken. Da kamen n\"{a}mlich ein
\emph{Punktsack} $\mathfrak P$ und ein \emph{Geradensack} $\mathfrak G$ zum
Vorschein, gemeinsam mit einer \emph{Inzidenz} genannten Teilmenge von
$\mathfrak P \times \mathfrak G$. Dann wurden drei Axiome pr\"{a}sentiert, und
fertig war die Definition einer projektiven Ebene! Zur Abrundung gab es noch
drei Modelle: Die projektiv abgeschlossene Anschauungsebene, die
Sieben-Punkte-Ebene von Fano und das B\"{u}ndelmodell der gew\"{o}hnlichen projektiven
Ebene, in dem zur allgemeinen Verwirrung \"{u}bliche Geraden als
{\glqq}Punkte{\grqq} und \"{u}bliche Ebenen als {\glqq}Geraden{\grqq} zu bezeichnen
waren. Kurz gesagt: Es versprach spannend zu werden. Und es wurde spannend!
\par
In den folgenden Jahren h\"{o}rte ich bei Brauner Vorlesungen \"{u}ber
\emph{Differentialgeometrie}, \emph{H\"{o}here Differentialgeometrie},
\emph{Liniengeometrie} und \emph{Abbildungsverfahren der konstruktiven
Geometrie}. Als junger Assistent begleitete ich ihn auch in die Vorlesungen
\"{u}ber \emph{Darstellende Geometrie\/} f\"{u}r Studierende der Architektur, des
Bauingenieurwesens und der Geod\"{a}sie.
\begin{figure}[!t] \centering
    \includegraphics[width=\textwidth]{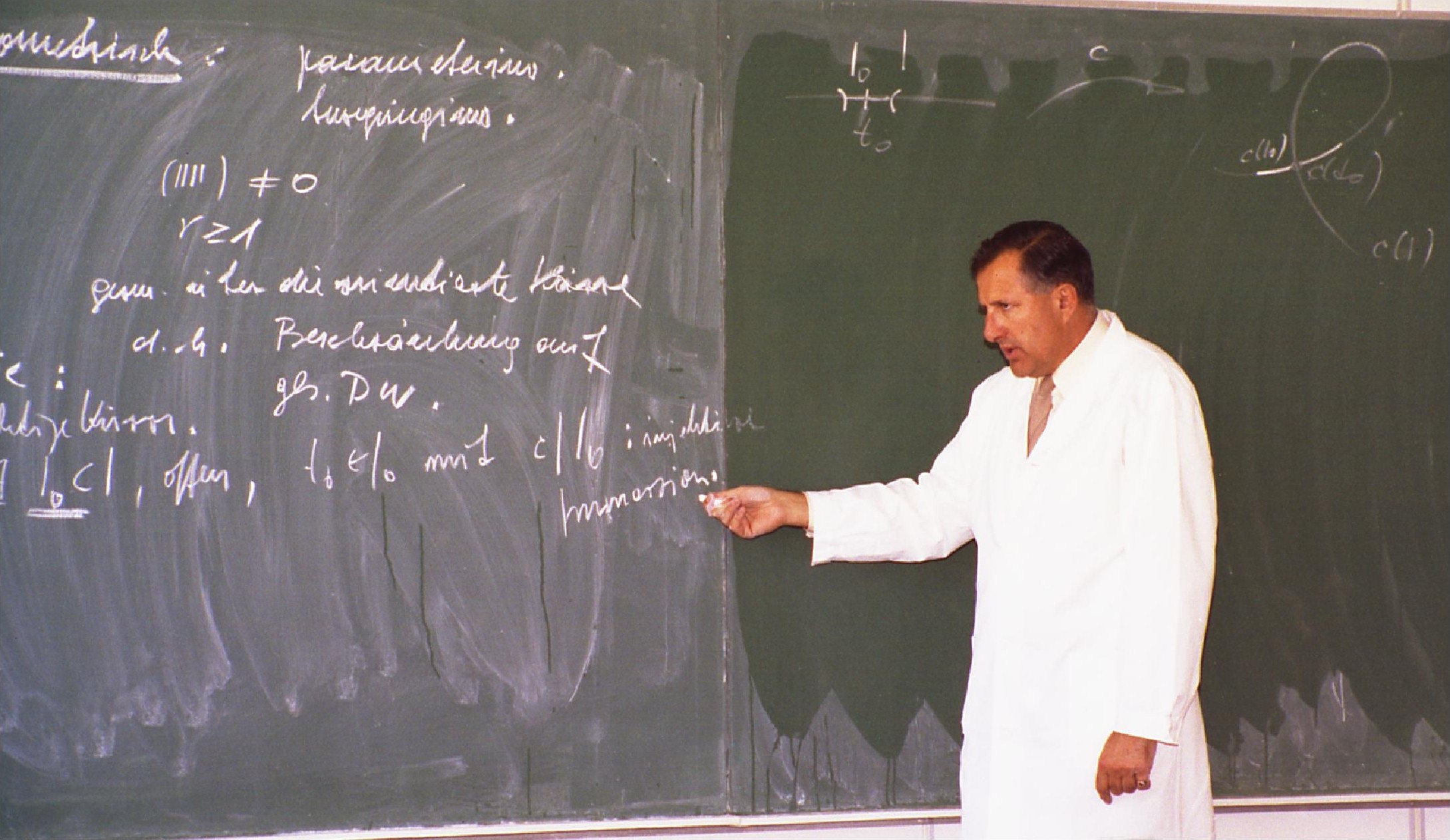}
    \caption{Vorlesung \"{u}ber Differentialgeometrie 1983}\label{bild4}
\end{figure}
\par
\begin{figure}[!t] \centering
    \includegraphics[width=0.9\textwidth]{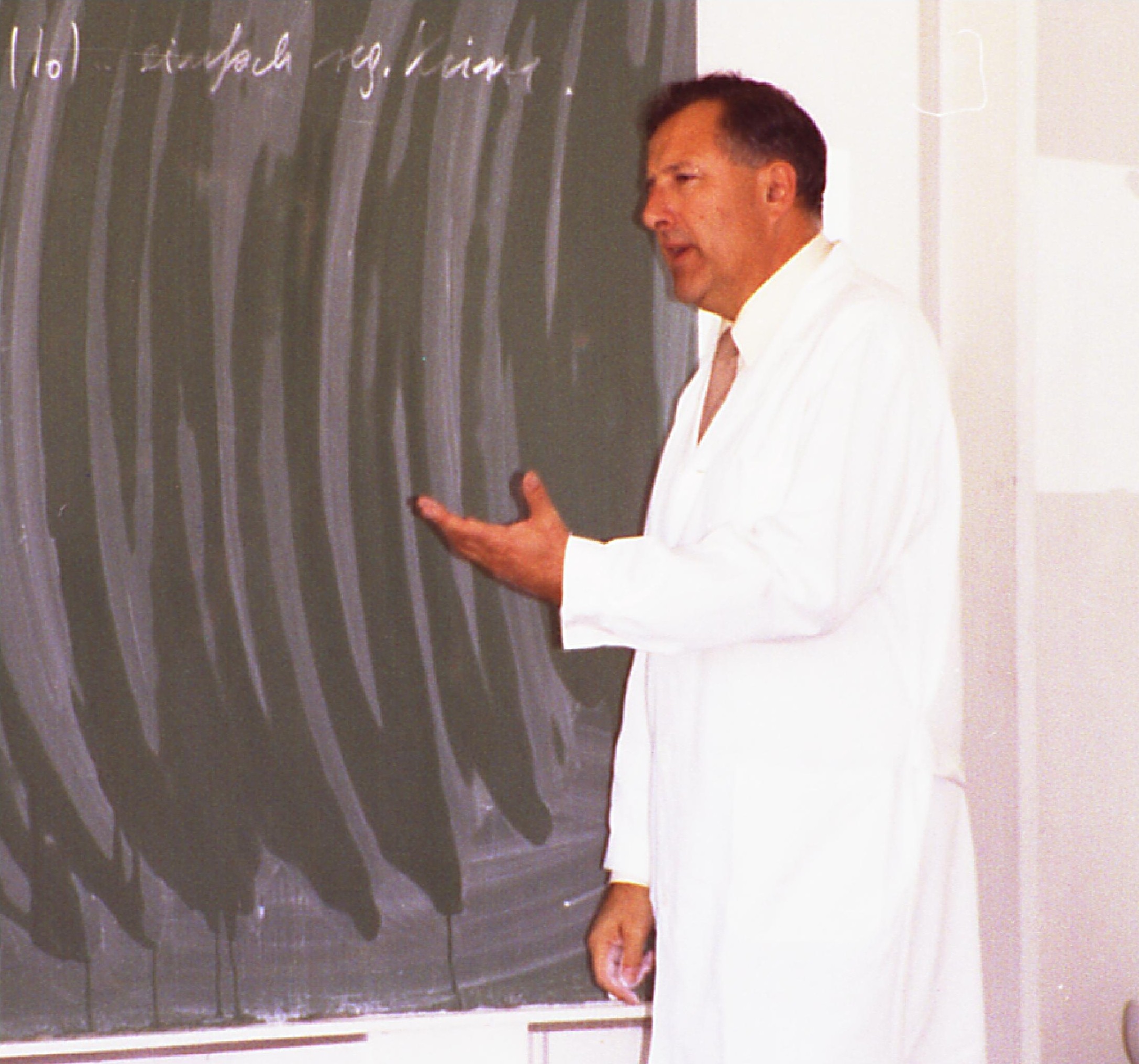}
    \caption{Vorlesung \"{u}ber Differentialgeometrie 1983}\label{bild5}
\end{figure}

Leider gibt es nur ganz wenige Bilder aus Brauners Lehrveranstaltungen. Das
Foto in Abbildung \ref{bild3} habe ich im Sommersemester 1982 in einer
Vorlesung \"{u}ber Differentialgeometrie aufgenommen, die an der Technischen
Universit\"{a}t Wien stattfand. Zu sehen ist Brauner, wie er gerade den Hauptsatz
der Hyperfl\"{a}chentheorie beweist: \emph{Zwei Immersionen mit derselben ersten
und zweiten Fundamentalform sind bewegungsgleich.} Im Wintersemester 1983
entstanden ebenfalls an der Technischen Universit\"{a}t Wien in einer der ersten
Vorlesungen \"{u}ber Differentialgeometrie jene beiden Aufnahmen, die in den
Abbildungen \ref{bild4} und \ref{bild5} zu sehen sind: Brauner erkl\"{a}rt hier,
den Blick ins Unendliche gerichtet, was unter einer \emph{geometrischen Aussage
\"{u}ber eine Kurve\/} zu verstehen sei.
\par
Brauners Vorlesungen \"{u}ber \emph{Differentialgeometrie} sind mir in bester
Erinnerung. Sie fanden w\"{a}hrend meiner Studienzeit nicht an der Technischen
Hochschule Wien, sondern an der Universit\"{a}t Wien in den R\"{a}umen des
Priesterseminars statt. Brauner setzte von Anfang an voraus, dass man Analysis
und Lineare Algebra \emph{schon gelernt hatte}. Das f\"{u}hrte dazu, dass bereits
in der zweiten Vorlesung deutlich weniger H\"{o}rerinnen und H\"{o}rer waren, als in
der ersten. So fand wenigstens ab diesem Zeitpunkt jeder einen Sitzplatz. Das
war auch gut so. Da es n\"{a}mlich kein Skriptum gab, mussten wir auf den kleinen,
an den H\"{o}rsaalst\"{u}hlen angebrachten Klapptischchen all das mitschreiben, was
Brauner rasant auf der Tafel notierte.
\par
Die Zielsetzung f\"{u}r diese Vorlesungen kann auch heute noch in der Einleitung
seines Lehrbuches der Differentialgeometrie nachgelesen werden. Dort schreibt
Brauner:
\begin{center}
\begin{minipage}[t]{0.9\textwidth}
{\glqq}Differentialgeometrie ist meines Erachtens ein Gebiet, das sich wegen
zahlreicher Querverbindungen zu anderen mathematischen Disziplinen und seiner
Bedeutung etwa f\"{u}r die theoretische Physik besonders gut als Vorlesung f\"{u}r den
zweiten Studienabschnitt einer Mathematikerausbildung eignet.{\grqq}
\end{minipage}
\end{center}
\par
Wie ich zuvor schon andeutete, war Brauner in seinen Vorlesungen sehr z\"{u}gig
unterwegs. Langatmige Motivationen oder ausgedehnte Wiederholungen des Stoffes
waren ihm fremd. Dennoch war es einfach faszinierend und vor allem lohnend,
seine Vorlesungen zu besuchen. Mathematisch pr\"{a}zise Formulierungen und
glasklare Definitionen, gepaart mit anschaulich-geometrischen Erkl\"{a}rungen und
zahlreichen Handskizzen, bildeten die Basis seiner Vorlesungen. Wer regelm\"{a}{\ss}ig
seine Veranstaltungen besuchte, wusste immer, worum es gerade ging.
\par
Eine seiner besonderen Eigenarten war es, beim Fenster hinaus blickend zu
unterrichten. In solchen Augenblicken wussten wir Studenten: Jetzt ist er voll
bei der Sache; nichts und niemand kann ihn aufhalten. Aber gelegentlich hielt
Brauner von sich aus pl\"{o}tzlich inne, dachte wortlos nach, sch\"{u}ttelte manchmal
auch den Kopf, schwieg nochmals f\"{u}r einige Sekunden, um dann im gewohnten Tempo
weiterzumachen.
\par
In seinen Lehrveranstaltungen konnte Brauner begeistern und mitrei{\ss}en. Einer
meiner Studienkollegen wollte im Anschluss an ein Seminar im Studienjahr
1975/76 zum Thema \emph{Nichtdesarguessche Projektive Ebenen\/} in seiner
Hausarbeit unbedingt das Problem der Existenz oder Nichtexistenz einer
projektiven Ebene der Ordnung $10$ l\"{o}sen. Brauner, der um die extreme
Schwierigkeit der Fragestellung wusste, hat ihm mit Recht ein anderes Thema
vorgeschlagen. Das genannte Problem wurde \"{u}brigens von Lam, Thiel und Swiercz
erst 1989 unter Einsatz des Computers gel\"{o}st. Wir wissen seither, dass es keine
solche Ebene gibt.
\par
Brauners Vortragsstil war in jeder Hinsicht brillant. Mit wenigen, treffenden
Worten das Richtige zu sagen, das war eine seiner St\"{a}rken. Er sprach laut,
deutlich und in ganzen S\"{a}tzen, die in vielen F\"{a}llen druckreif waren. Ich
erinnere mich an einen Artikel in einer Studentenzeitung aus den 1980er Jahren.
Dort wurde Brauner als der {\glqq}ungekr\"{o}nte Meister des Schachtelsatzes{\grqq}
bezeichnet. Dem habe ich nichts hinzuzuf\"{u}gen.
\par
Gelegentlich streute Brauner in seinen Unterricht aber auch launische
Bemerkungen ein. So erkl\"{a}rte er in einer Vorlesung \"{u}ber Differentialgeometrie
die kovariante Ableitung auf einer Fl\"{a}che mit Hilfe von {\glqq}auf einer Fl\"{a}che
lebenden K\"{a}fern{\grqq} und bemerkte dabei verschmitzt:
\begin{center}
%\begin{minipage}[t]{0.9\textwidth}
    {\glqq}Nur differenzieren sollten die K\"{a}fer schon k\"{o}nnen.{\grqq}
%\end{minipage}
\end{center}
In einer Vorlesung direkt vor den Osterferien schrieb er zum Abschluss
\begin{center}
    \textsf{Frohe O*}
\end{center}
auf die Tafel, um dann wortlos schmunzelnd den Raum zu verlassen. In seinen
Vorlesungen f\"{u}r Ingenieurstudenten betonte er zur Illustration eines r\"{a}umlichen
Rechtssystems immer wieder nachdr\"{u}cklich:
\begin{center}
\begin{minipage}[t]{0.9\textwidth}
    {\glqq}Die $z$-Achse weist nach oben, die $y$-Achse nach weist nach rechts, und
    die $x$-Achse sticht Sie in den Bauch.{\grqq}
\end{minipage}
\end{center}
Brauner \"{u}bersetzte in einer Vorlesung aus projektiver Geometrie das Wort
\emph{oskulieren} korrekt als \emph{k\"{u}ssen} und meinte danach nur trocken:
\begin{center}
%\begin{minipage}[t]{0.9\textwidth}
    {\glqq}Was \emph{hyperoskulieren} bedeutet, m\"{u}ssen Sie selbst
    herausfinden.{\grqq}
%\end{minipage}
\end{center}
\par
Viele Inhalte der Vorlesungen von Heinrich Brauner k\"{o}nnen auch heute noch in
den sechs von ihm verfassten B\"{u}chern nachgelesen werden. Sie behandeln die
Themen \emph{Geometrie Projektiver R\"{a}ume} (2 B\"{a}nde), \emph{Baugeometrie} (2
B\"{a}nde, gemeinsam mit Walter Kickinger), \emph{Differentialgeometrie} und
\emph{Konstruktive Geometrie}. Seinen au{\ss}ergew\"{o}hnlichen Vortragsstil k\"{o}nnen sie
aus meiner Sicht leider nicht vermitteln.

\section{Der Forscher Heinrich Brauner}

\begin{figure}[!t]
\centering
    \includegraphics[width=0.90\textwidth]{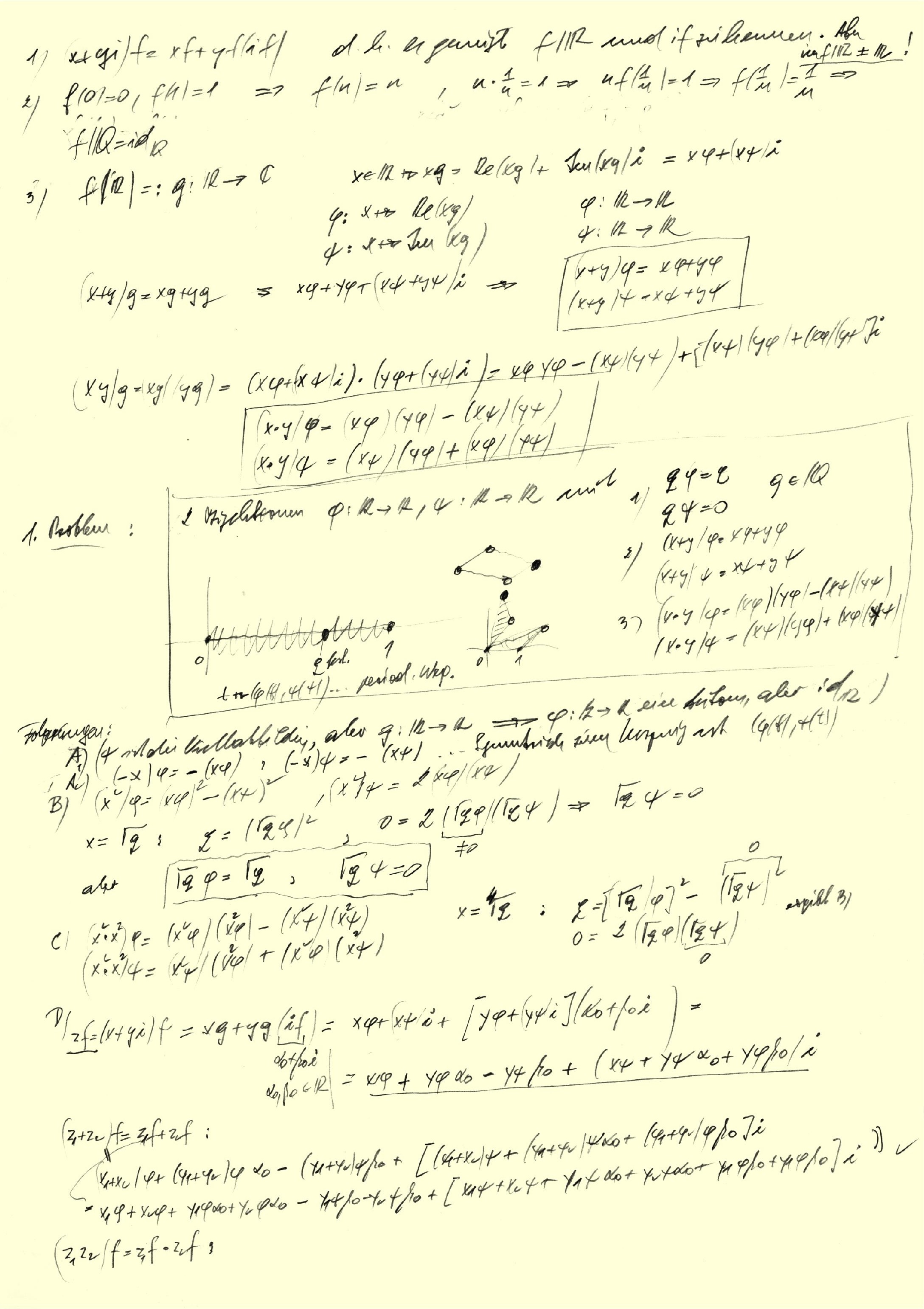}
    \caption{Ein von Brauner verfasstes Manuskript}\label{bild6}
\end{figure}

Brauner arbeitete an seinen Artikeln und B\"{u}chern weitgehend alleine und
vorzugsweise daheim. Er gab aber die von seiner Sekret\"{a}rin mit der
Schreibmaschine ausgearbeiteten Manuskripte immer uns Assistenten zum
Durchlesen, Kommentieren und Korrigieren. Aus diesem Grund sind nur sehr wenige
handgeschriebene Aufzeichnungen von Brauner vorhanden. Abbildung \ref{bild6}
zeigt ein Manuskript aus dem Jahre 1986, in dem er sich mit den Derivationen
des komplexen Zahlk\"{o}rpers besch\"{a}ftigte. Er hat dar\"{u}ber aber nichts publiziert.
\par
Umgekehrt nahm sich aber Brauner auch immer sehr viel Zeit, um die Artikel
seiner Mitarbeiter gewissenhaft zu studieren und zu verbessern. So manches
meiner Manuskripte war kaum mehr zu erkennen, nachdem es Brauner gelesen und --
wie immer mit Bleistift -- seine Anmerkungen angebracht hatte. Seine Kritik
bezog sich dabei prim\"{a}r auf den mathematischen Inhalt, wo er bei anderen
dieselben strengen Ma{\ss}st\"{a}be ansetzte wie bei sich selbst. Er markierte aber
prinzipiell alles, was ihm falsch erschien. Oft formulierte er seine
Bemerkungen zus\"{a}tzlich sehr pointiert, aber niemals unh\"{o}f"|lich, im
pers\"{o}nlichen Gespr\"{a}ch. So war etwa sein trockener Kommentar, nachdem er das
erste Kapitel meiner Dissertation gelesen hatte: {\glqq}Herr Havlicek, ihre
Beistrichsetzung m\"{o}chte ich nicht haben!{\grqq}
\begin{figure}[h!]
\centering
    \includegraphics[width=10.5cm]{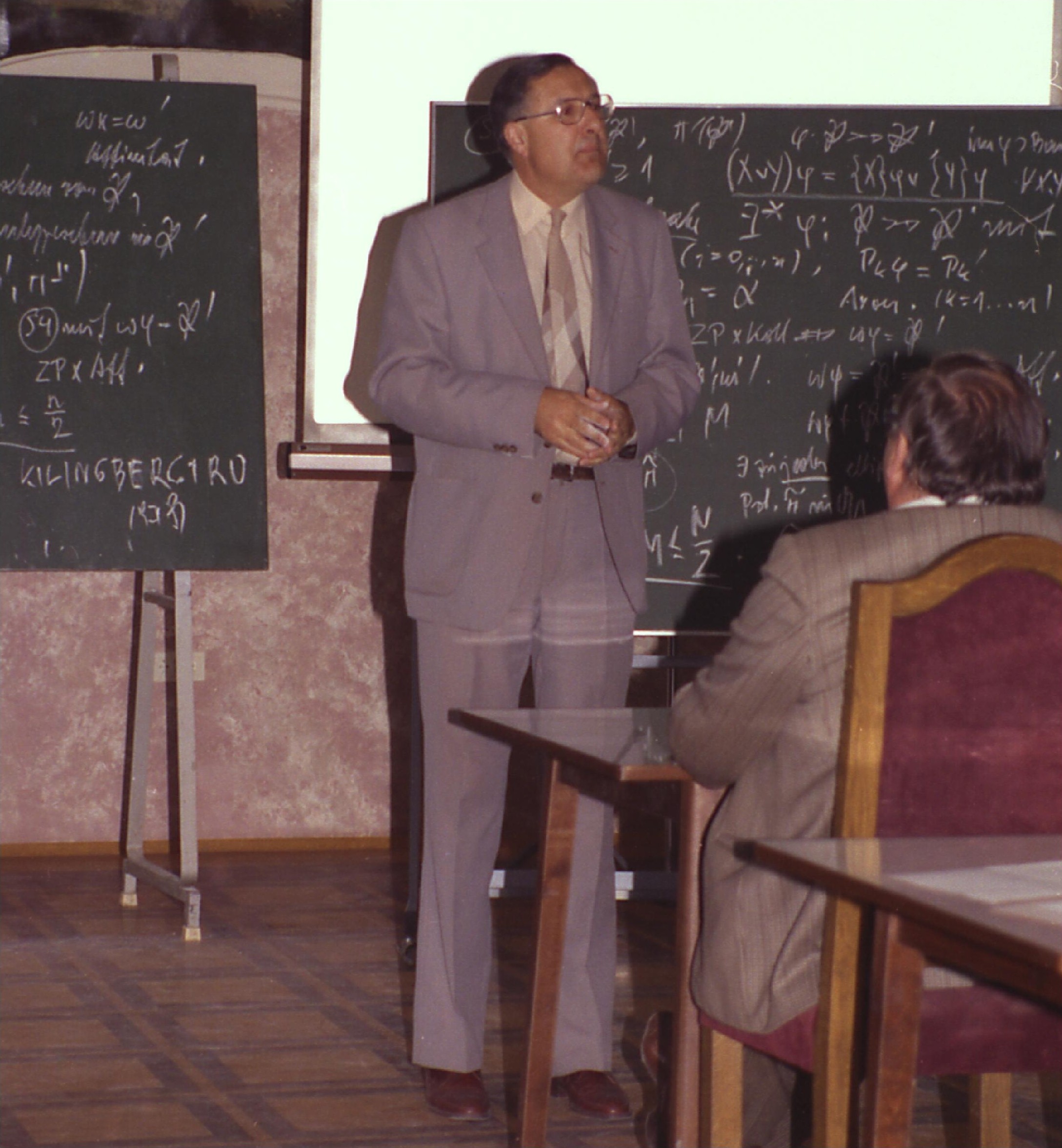}
    \caption{Vortrag im Stift Rein 1983}\label{bild7}
\end{figure}
\par
Brauner legte immer allergr\"{o}{\ss}ten Wert auf wissenschaftliche Gespr\"{a}che mit
seinen Mitarbeitern. So knapp konnte seine Zeit gar nicht bemessen sein, dass
er daf\"{u}r nicht ein paar Minuten er\"{u}brigen konnte. Und so manche Unterredung hat
dann deutlich l\"{a}nger gedauert, als urspr\"{u}nglich geplant war. So sprachen wir
einmal sicher f\"{u}r mehr als eine halbe Stunde -- auch wenn es unglaubw\"{u}rdig
klingen mag -- \"{u}ber die \emph{leere Menge}.
\par
\begin{figure}[h!]
\centering
    \includegraphics[width=\textwidth]{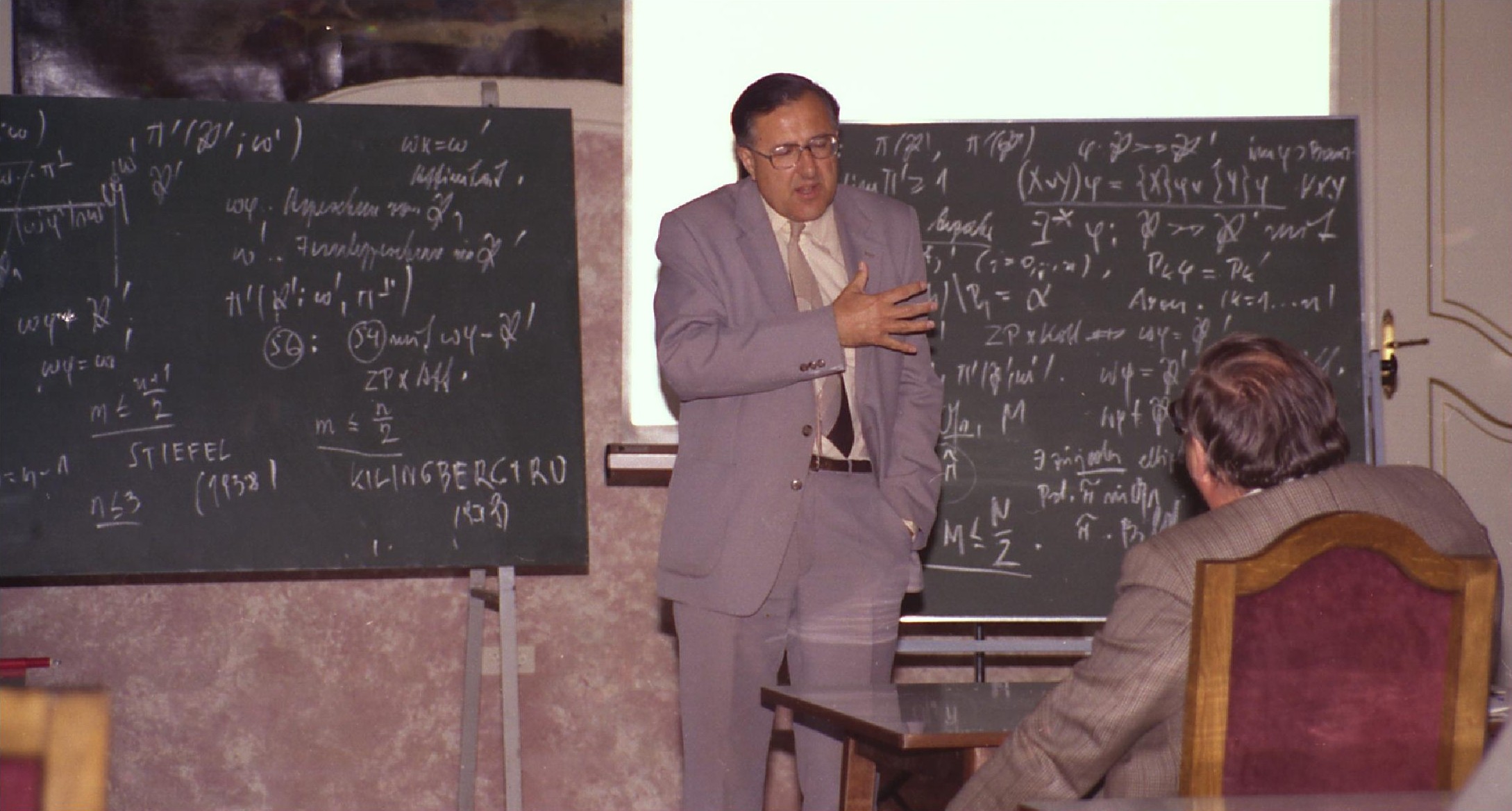}
    \caption{Vortrag im Stift Rein 1983}\label{bild8}
\end{figure}

Selbstverst\"{a}ndlich hat Brauner seine Forschungsergebnisse auf Tagungen
pr\"{a}\-sen\-tiert. Alles das, was ich zuvor \"{u}ber seine Vorlesungen angemerkt
habe, trifft auch auf den Stil seiner Vortr\"{a}ge zu. Die Abbildungen \ref{bild7}
und \ref{bild8} zeigen ihn beim \emph{Zweiten \"{O}sterreichischen
Geometrie-Kolloquium}, welches im Mai 1983 im Stift Rein stattfand. Er sprach
damals \"{u}ber den \emph{Satz von Pohlke im $n$-dimensionalen euklidischen Raum}.
\par
Brauners sehr breit gestreutes wissenschaftliches Werk hat Walter Wunderlich in
seinem Nachruf \cite{3} ausf\"{u}hrlich gew\"{u}rdigt. Im Anhang~1 zu diesem Artikel
ist ein Schriftenverzeichnis so wiedergegeben, wie es Brauner selbst gef\"{u}hrt
hat. Es umfasst einundneunzig Arbeiten.

\section{Schlussbemerkungen}

Es g\"{a}be noch viel zu berichten, etwa \"{u}ber die zwanzig von Brauner betreuten
Dissertationen (vgl.\ dazu den Anhang~2) und die wohl mehr als einhundert von
ihm vergebenen Haus- und Diplomarbeiten, \"{u}ber welche es allerdings keine
vollst\"{a}ndigen schriftlichen Aufzeichnungen geben d\"{u}rfte.
\par
Neben seinen Aktivit\"{a}ten in der universit\"{a}ren Lehre und Forschung galt sein
gro{\ss}es Engagement insbesondere dem Unterrichtsfach \emph{Darstellende
Geometrie}, und zwar in inhaltlicher, didaktischer und fachpolitischer
Hinsicht. Auch dieser Aspekt muss hier leider ausgeklammert bleiben.
\par
Hingegen m\"{o}chte ich zum Abschluss die gro{\ss}e Hilfsbereitschaft Brauners in
gro{\ss}en wie in kleinen Dingen erw\"{a}hnen. Dazu sei eine der Situationen
geschildert, in denen mir Brauner entscheidend geholfen hat:
\par
Zum Ende meines Studiums hatte ich Brauner als Pr\"{u}fer f\"{u}r die m\"{u}ndliche
Lehramtspr\"{u}fung aus Darstellender Geometrie gew\"{a}hlt. Aber genau eine Woche vor
dieser Pr\"{u}fung h\"{a}tte ich eine Truppen\"{u}bung beim \"{o}sterreichischen Bundesheer
absolvieren m\"{u}ssen. Ich stellte einen Antrag auf Aufschub der Einberufung und
gab als Begr\"{u}ndung an, dass ich mich in Ruhe auf meine Abschlusspr\"{u}fung
vorbereiten m\"{u}sste. Mein Antrag wurde abgelehnt, da keine Terminkollision
vorlag. Einer meiner Studienkollegen, dem Brauner in anderem Zusammenhang
erfolgreich geholfen hatte, empfahl mir, in Brauners Sprechstunde zu gehen. Ich
folgte seinem Rat und schilderte Brauner mein Problem. Dieser setzte einfach
eine {\glqq}Vorpr\"{u}fung{\grqq} mitten im Zeitraum der Truppen\"{u}bung an und
verfasste eine Best\"{a}tigung dar\"{u}ber. Mit diesem Schreiben legte ich erfolgreich
Berufung gegen die Entscheidung der ersten Instanz ein.
\par
Damit bin ich am Ende meiner Ausf\"{u}hrungen angelangt. Im Mittelpunkt standen
meine ganz pers\"{o}nlichen Erinnerungen an den Menschen Heinrich Brauner, seine
Eigenschaften, Wesensmerkmale und Leistungen. Wer mehr erfahren m\"{o}chte, dem
seien neben seinem wissenschaftlichen Werk (vgl.\ das Schriftenverzeichnis im
Anhang~1) auch der Artikel \cite{1}, die Ausarbeitung seiner Wiener
Antrittsvorlesung vom 28.~J\"{a}nner 1970 \cite{2} und der schon erw\"{a}hnte Nachruf
\cite{3} w\"{a}rmstens ans Herz gelegt.

%% Ab hier fallweise alte Rechtschreibung, zB Fallinie.

\section*{Anhang 1: Publikationen von H.~Brauner}

\begin{publikationen}{99}

\item
Orthogonalsysteme von Riemannschen Hyperfl\"{a}chen der Klasse 1. \emph{Anz.\
\"{O}ster.\ Akad.\ Wiss.\ Math.-Nat.\ Kl.} \textbf{88} (1951). 29--36.

\item
Kongruente Verlagerung kollinearer R\"{a}ume in axiale Lage. \emph{Monatsh.\ Math.}
\textbf{57} (1953). 75--87.

\item
Kongruente Verlagerung kollinearer R\"{a}ume in halbxiale Lage. \emph{Monatsh.\
Math.} \textbf{58} (1954). 13--26.

\item
Quadriken als Bewegfl\"{a}chen. \emph{Monatsh.\ Math.} \textbf{59} (1955). 45--63.

\item
Erzeugung eines gleichseitigen hyperbolischen Paraboloides durch Bewegung einer
gleichseitigen Hyperbel. \emph{Arch.\ Math.\ (Basel)\/} \textbf{6} (1955).
330--334.

\item
Geod\"{a}tische Fallinien einer Gel\"{a}ndefl\"{a}che. \emph{Anz.\ \"{O}ster.\ Akad.\ Wiss.\
Math.-Nat.\ Kl.} \textbf{92} (1955). 171--175.

\item
\"{U}ber die Projektion mittels der Sehnen einer Raumkurve 3.~Ordnung.
\emph{Monatsh.\ Math.} \textbf{59} (1955), 258--273.

\item
\"{U}ber die \"{a}hnlichen und sich \"{a}hnlich projizierenden Kegelschnitte auf Quadriken.
\emph{Arch.\ Math.\ (Basel)\/} \textbf{7} (1956), 78--86.

\item
Konstruktive Durchf\"{u}hrung der durch die Sehnen einer Raumkurve 3.~Ordnung
vermittelten Abbildung des Raumes auf eine Ebene. \emph{Monatsh.\
Math.}\textbf{60} (1956), 231--248.

\item
Die automorphen involutorischen Korrelationen koaxialer projektiver
Schraubungen (mit Rudolf Bereis). \emph{\"{O}sterreich.\ Akad.\ Wiss.\ Math.-Nat.\
Kl.\ S.-B.\ II.} \textbf{165} (1956), 327--355.

\item
\"{U}ber Mannigfaltigkeiten von Strahlen mit kongruenten Netzrissen. \emph{Arch.\
Math.\ (Basel)\/} \textbf{7} (1957), 406--416.

\item
\"{U}ber koaxiale euklidische Schraubungen (mit Rudolf Bereis). \emph{Monatsh.\
Math.} \textbf{61} (1957), 225--245.

\item
Schraubung und Netzprojektion. \emph{Elem.\ Math.} \textbf{12} (1957). 33--41.
\item
Eine Verallgemeinerung der Zyklographie. \emph{Arch.\ Math.\ (Basel)\/}
\textbf{9} (1958), 470--480.

\item
\"{U}ber die durch einen quadratischen Komplex der Charakteristik (11)(112)
vermittelte Projektion I. \emph{Monatsh.\ Math.} \textbf{62} (1958), 119--131.

\item
\"{U}ber die durch einen quadratischen Komplex der Charakteristik (11)(112)
vermittelte Projektion II. \emph{Monatsh.\ Math.} \textbf{62} (1958), 132--145.

\item
Bestimmung einer Strahlfl\"{a}che aus ihren sph\"{a}rischen Bildern. \emph{Anz.\
\"{O}ster.\ Akad.\ Wiss.\ Math.-Nat.\ Kl.} \textbf{95} (1958). 103--107.

\item
\"{U}ber Strahlfl\"{a}chen von konstantem Drall. \emph{Monatsh.\ Math.} \textbf{63}
(1959), 101--111.

\item
Die dualen Gegenst\"{u}cke zu fl\"{a}chentheoretischen S\"{a}tzen von O.~Bonnet und
E.~Beltrami. \emph{Anz.\ \"{O}ster.\ Akad.\ Wiss.\ Math.-Nat.\ Kl.} \textbf{96}
(1959), 194--200.

\item
Eine Verallgemeinerung des Problems der Ces\`{a}rokurven. \emph{Math.\ Ann.}
\textbf{138} (1959), 27--41.

\item
Beitr\"{a}ge zur Theorie des mit einer euklidischen Schraubung verkn\"{u}pften
kubischen Nullsystems (mit Rudolf Bereis). \emph{Math.\ Nachr.} \textbf{20}
(1959), 239--258.

\item
Die Strahlfl\"{a}che 3. Grades mit konstantem Drall. \emph{Monatsh.\ Math.}
\textbf{64} (1960), 101--109.

\item
Erweiterung des Begriffes Drall auf Mongesche Fl\"{a}chen. \emph{Anz.\ \"{O}ster.\
Akad.\ Wiss.\ Math.-Nat.\ Kl.} \textbf{97} (1960). 139--144.

\item
Die konstant gedrallte Netzfl\"{a}che 4.~Grades. \emph{Monatsh.\ Math.} \textbf{65}
(1961), 53--73.

\item
Eine einheitliche Erzeugung konstant gedrallter Strahlfl\"{a}chen. \emph{Monatsh.\
Math.} \textbf{65} (1961), 301--314.

\item
Die verallgemeinerten B\"{o}schungsfl\"{a}chen. \emph{Math.\ Ann.} \textbf{143} (1961),
431--439.

\item
Die Affinnormalen der Tangentialschnitte einer Fl\"{a}che. \emph{Anz.\ \"{O}ster.\
Akad.\ Wiss.\ Math.-Nat.\ Kl.} \textbf{99} (1962). 9--14.

\item
Eine Scherungsinvariante der Strahlfl\"{a}chen. \emph{Monatsh.\ Math.} \textbf{66}
(1962), 105--109.

\item
Die konstant gedrallten windschiefen Fl\"{a}chen 4.~Grades mit reduzibler
Fernkurve. \emph{Math.\ Z.} \textbf{82} (1963), 420--433.

\item
Die windschiefen Fl\"{a}chen konstanter konischer Kr\"{u}mmung. \emph{Math.\ Ann.}
\textbf{152} (1963), 257--270.

\item
Geometrie auf der Cayleyschen Fl\"{a}che. \emph{\"{O}sterreich.\ Akad.\ Wiss.\
Math.-Natur.\ Kl.\ S.-B.\ II\/} \textbf{173} (1964), 93--128.

\item
Kreisgeometrie in der isotropen Ebene. \emph{Monatsh.\ Math.} \textbf{69}
(1965), 105--128.

\item
Die quadratischen Strahlkomplexe der Charakteristik (321). \emph{Math.\ Z.}
\textbf{88} (1965), 320--357.

\item
Geometrie des zweifach isotropen Raumes. I. Bewegungen und kugeltreue
Transformationen. \emph{J.\ Reine Angew.\ Math.} \textbf{224} (1966), 118--146.

\item
Die Fl\"{a}chen mit einem kinematischen Netz aus Schmieglinien (mit Hermann
Schaal). \emph{Arch.\ Math.\ (Basel)\/} \textbf{18} (1967), 91--99.

\item
Geometrie des zweifach isotropen Raumes II. Differentialgeometrie der Kurven
und windschiefen Fl\"{a}chen. \emph{J.\ Reine Angew.\ Math.} \textbf{226} (1967),
132--158.

\item
Die algebraischen windschiefen Gesimsfl\"{a}chen. \emph{Monatsh.\ Math.}
\textbf{71} (1967), 300--318.

\item
Geometrie des zweifach isotropen Raumes III. Fl\"{a}chentheorie. \emph{J.\ Reine
Angew.\ Math.} \textbf{228} (1967), 38--70.

\item
\emph{Differentialgeometrie.} Universit\"{a}t Stuttgart, Stuttgart 1967.
vi+127~Seiten.

\item
Neuere Untersuchungen \"{u}ber windschiefe Fl\"{a}chen: Ein Bericht. \emph{Jber.\
Deutsch.\ Math.-Verein.} \textbf{70} (1967) Heft 2, Abt.\ 1, 61--85.

\item
\emph{Analytische Geometrie I.} Universit\"{a}t Stuttgart, Stuttgart 1967.
vi+114~Seiten.

\item
Die algebraischen windschiefen Fl\"{a}chen mit einer stetigen Schar ebener
Schattengrenzen. \emph{Math.\ Ann.} \textbf{176} (1968), 1--14.

\item
\emph{Analytische Geometrie II.} Universit\"{a}t Stuttgart, Stuttgart 1968.
iv+90~Seiten.

\item
\emph{Analytische Geometrie III.} Universit\"{a}t Stuttgart, Stuttgart 1968.
ii+223~Seiten.

\item
Die Fl\"{a}chen mit B\"{o}schungslinien als Fallinien. \emph{Monatsh.\
Math.}\textbf{72} (1968), 385--411.

\item
Die Fl\"{a}chen mit zwei Scharen konstant geb\"{o}schter Schmieglinien (mit Hermann
Schaal). \emph{Arch.\ Math.\ (Basel)\/} \textbf{20} (1969), 81--87.

\item
Die windschiefen Kegelschnittfl\"{a}chen. \emph{Math.\ Ann.} \textbf{183} (1969),
33--44.

\item
Die Fl\"{a}chen, welche stetige Scharen ebener geod\"{a}tischer Linien tragen.
\emph{Jber.\ Deutsch.\ Math.-Verein.} \textbf{71} (1969), Heft 3, Abt. 1,
160--166.

\item
\emph{Differentialgeometrie.} Universit\"{a}t Stuttgart, Stuttgart 1969.
ii+336~Seiten.

\item
\emph{Riemannsche Geometrie.} Universit\"{a}t Stuttgart, Stuttgart 1969.
136~Seiten.

\item
Gedanken \"{u}ber Geometrie. \emph{Antrittsvorlesungen der Technischen Hochschule
Wien\/} \textbf{12}. Verlag der Technischen Hochschule Wien, Wien 1970.
11~Seiten.

\item
Differentialgeometrie ebener Kurven.\emph{Wiss.\ Nachrichten\/} \textbf{26}
(1971), 21--24.

\item
Eine geometrische Kennzeichnung linearer Abbildungen. \emph{Monatsh.\ Math.}
\textbf{77} (1973), 10--20.

\item
Abbildungsmethoden der konstruktiven Geometrie. 7.\ Steierm\"{a}rkisches
Mathematisches Symposium (Graz, 1975), \emph{Ber.\ Math.-Statist.\ Sektion,
Forschungszentrum Graz\/} Nr.\ \textbf{38} (1975). 11~Seiten.
\item
\emph{Geometrie projektiver R\"{a}ume I.} Bibliographisches Institut,
Mannheim-Wien-Z\"{u}rich 1976. x+225~Seiten. (ISBN~10: 3-411-01512-8).

\item
\emph{Geometrie projektiver R\"{a}ume II.} Bibliographisches Institut,
Mannheim-Wien-Z\"{u}rich 1976. viii+250~Seiten. (ISBN~10: 3-411-01513-6).

\item
\emph{Baugeometrie I\/} (mit Walter Kickinger). 1.\ Auf"|lage.
Wiesbaden-Berlin, Bauverlag 1977. 88~Seiten. (ISBN~10: 3-762-50825-9).

\item
\"{U}ber schmieglinientreue Isometrien. \emph{\"{O}sterreich.\ Akad.\ Wiss.\
Math.-Natur.\ Kl.\ Sitzungsber.\ II\/} \textbf{188} (1979), no.\ 1--3, 15--21.

\item
Die erzeugendentreuen konformen Abbildungen aus Regelfl\"{a}chen. \emph{Arch.\
Math.\ (Basel)\/} \textbf{33} (1979/80), no.\ 5, 470--477.

\item
\emph{Geometrija u Graditeljstvu\/} (mit Walter Kickinger), \v{S}kolska knjiga,
Zagreb 1980. 156~Seiten. (\"{U}bersetzung von \emph{Baugeometrie I}, in kroatischer
Sprache).

\item
Abbildungen aus Regelfl\"{a}chen. 12.\ Steierm\"{a}rkisches Mathematisches Symposium
(Graz, 1980), \emph{Ber.\ Math.-Statist.\ Sektion, Forschungszentrum Graz\/}
Nr.\ \textbf{140} (1980). 14~Seiten.

\item
\emph{Differentialgeometrie.} Friedr.\ Vieweg \& Sohn, Braunschweig 1981.
xvii+424~Seiten. (ISBN~10: 3-528-03809-8).

\item
Gedanken zum Unterricht in Darstellender Geometrie. \emph{\"{O}MG Didaktik-Reihe\/}
\textbf{6} (1981). 76~Seiten.

\item
Darstellende Geometrie im Schulunterricht. \emph{Mathematikunterr.} \textbf{27}
(3), (1981), 5--68.

\item
Die fl\"{a}chentreuen Abbildungen aus Regelfl\"{a}chen, bei denen die Erzeugenden
geradlinig bleiben. \emph{Arch.\ Math.\ (Basel)\/} \textbf{38} (1982), no.\ 2,
102--105.

\item
\emph{Baugeometrie II\/} (mit Walter Kickinger). Wiesbaden-Berlin, Bauverlag
1982. 89~Seiten. (ISBN~10: 3-7625-0927-1).

\item
Gebaute Geometrie. Beispiele aus dem Bauwesen f\"{u}r den Schulunterricht der
Darstellenden Geometrie (mit Walter Kickinger). \emph{Mathematikunterr.}
\textbf{28} (2) (1982), 5--28.

\item
Zur Theorie linearer Abbildungen. \emph{Abh.\ Math.\ Sem.\ Univ.\ Hamburg\/}
\textbf{53} (1983), 154--169.

\item
Die windschiefen Fl\"{a}chen mit B\"{o}schungsschmieglinien. \emph{Anz. \"{O}sterreich.\
Akad.\ Wiss.\ Math.-Natur.\ Kl.} \textbf{121} (1984), 125--127 (1985).

\item
Zur theoretischen Begr\"{u}ndung der Darstellenden Geometrie. \emph{Ber.\
Math.-Statist.\ Sektion, Forschungszentrum Graz\/} Nr.\ \textbf{227} (1985).
2~Seiten.

\item
Zur Methodik der Darstellenden Geometrie I. Die konstruktive Behandlung der
Ebene. \emph{Informationsbl\"{a}tter Darstellende Geometrie\/} (Univ.\ Innsbruck)
\textbf{4} (1), (1985), 11--17.

\item
Die erzeugendentreuen geod\"{a}tischen Abbildungen aus Regelfl\"{a}chen.
\emph{Monatsh.\ Math.} \textbf{99} (1985), no.\ 2, 85--103.

\item
Zur Methodik der Darstellenden Geometrie II. Der Anfangsunterricht.
\emph{Informationsbl\"{a}tter Darstellende Geometrie\/} (Univ.\ Innsbruck)
\textbf{4} (2), (1985), 15--24.

\item
Lineare Abbildungen aus euklidischen R\"{a}umen. \emph{Beitr\"{a}ge Algebra Geom.}
\textbf{21} (1986), 5--26.

\item
Die verallgemeinerten B\"{o}schungsfl\"{a}chen mit B\"{o}schungsschmieglinien.
\emph{\"{O}sterreich.\ Akad.\ Wiss.\ Math.-Natur.\ Kl.\ Sitzungsber.\ II\/}
\textbf{194} (1985), no.\ 1--3, 55--61.

\item
Zur Methodik der Darstellenden Geometrie III. L\"{o}sung stereometrischer Aufgaben
mit Hilfe von Normalprojektionen. \emph{Informationsbl\"{a}tter Darstellende
Geometrie\/} (Univ.\ Innsbruck) \textbf{5} (1), (1986), 7--13.

\item
\emph{Lehrbuch der konstruktiven Geometrie.} Wien, Springer, 1986. 384~Seiten.
(ISBN~10: 3-211-81833-2).

\item
Zur Methodik der darstellenden Geometrie IV. Parallelri{\ss} einer Ellipse.
\emph{Informationsbl\"{a}tter Darstellende Geometrie\/} (Univ.\ Innsbruck)
\textbf{5} (2) (1986), 11--16.

\item
Eine Kennzeichnung der Minimalfl\"{a}chen von G.~Thomsen. \emph{Rad Jugoslav.\
Akad.\ Znan.\ Umjet.\/} no.\ \textbf{435}, (1988), 1--15.

\item
Zum Satz von K.~Pohlke in $n$-dimensionalen euklidischen R\"{a}umen.
\emph{\"{O}sterreich.\ Akad.\ Wiss.\ Math.-Natur.\ Kl.\ Sitzungsber.\ II\/}
\textbf{195} (1986), no.\ 8--10, 585--591.

\item
Zur Methodik der darstellenden Geometrie V. Methodische Miniaturen.
\emph{Informationsbl\"{a}tter Darstellende Geometrie\/} (Univ.\ Innsbruck)
\textbf{6} (1) (1987), 13--20.

\item
Darstellende Geometrie an der AHS -- ein Unterrichtsgegenstand im Wandel.
\emph{Informationsbl\"{a}tter Darstellende Geometrie\/} (Univ.\ Innsbruck)
\textbf{6} (1) (1987), 3--10.

\item
Zur Methodik der darstellenden Geometrie VI. Der Unterrichtsgegenstand
Darstellende Geometrie im Zeitalter des Computers. \emph{Informationsbl\"{a}tter
Darstellende Geometrie\/} (Univ.\ Innsbruck) \textbf{6} (2) (1987), 11--18.
\item
Die Drehfl\"{a}chen mit B\"{o}schungsschmieglinien. \emph{\"{O}sterreich.\ Akad.\ Wiss.\
Math.-Natur.\ Kl.\ Sitzungsber.\ II\/} \textbf{196} (1987), no.\ 4--7,
217--226.

\item
Zur Methodik der darstellenden Geometrie VII. Abbildungen im Unterricht der
darstellenden Geometrie, Teil 1. \emph{Informationsbl\"{a}tter Darstellende
Geometrie\/} (Univ.\ Innsbruck) \textbf{7} (1), (1988), 7--16.

\item
Eine Kennzeichnung der \"{A}hnlichkeiten affiner R\"{a}ume mit definiter
Orthogonalit\"{a}tsstruktur. \emph{Geom.\ Dedicata\/} \textbf{29} (1989), no.\ 1,
45--51.

\item
Zur Methodik der darstellenden Geometrie VIII. Abbildungen im Unterricht der
Darstellenden Geometrie, Teil 2, Abbildungsgleichungen zur Herstellung von
Rissen. \emph{Informationsbl\"{a}tter Darstellende Geometrie\/} (Univ.\ Innsbruck)
\textbf{7} (2), (1988), 13--22.

\item
\"{U}ber die von Kollineationen projektiver R\"{a}ume induzierten Geradenabbildungen.
\emph{\"{O}sterreich.\ Akad.\ Wiss.\ Math.-Natur.\ Kl.\ Sitzungsber.\ II\/}
\textbf{197} (1988), no.\ 4--7, 327--332.

\item
Zur Methodik der Darstellenden Geometrie IX. Abbildungen im Unterricht der
Darstellenden Geometrie, Teil 3. \emph{Informationsbl\"{a}tter Darstellende
Geometrie\/} (Univ.\ Innsbruck) \textbf{8} (1) (1989), 5--14.

\item
\emph{Baugeometrie I\/} (mit Walter Kickinger). 2.\ Auf"|lage. Wiesbaden
Berlin, Bauverlag 1989. 91~Seiten. (ISBN~10: 3-762-52690-7).

\item
Die Schraubfl\"{a}chen und die Spiralfl\"{a}chen mit B\"{o}schungsschmieglinien.
\emph{Glas.\ Mat.\ Ser.\ III\/} \textbf{25} (45) (1990), no.\ 1, 157--165.
\end{publikationen}

\section*{Anhang 2: Dissertationen\\(Betreuer und Erstgutachter H.~Brauner)}
\begin{publikationen}{99}

\item
Oswald Giering: \emph{Bestimmung von Eibereichen und Eik\"{o}rpern durch
Steiner-Sym\-metrisierungen.}
Technische Hochschule Stuttgart 1962.\\
Erschienen in: \emph{Sitz.-Ber.\ Bayer.\ Akad.\ Wiss.\ M\"{u}nchen\/} (1962),
225--253.

\item
Wolfgang Jenne: \emph{Eine nat\"{u}rliche Affingeometrie der Strahlfl\"{a}chen.}
Technische Hochschule Stuttgart, 1964.\\
Erschienen in: \emph{Math.\ Zeitschr.} \textbf{83} (1964), 214--237.

\item
Gerd Blind: \emph{Ebene Lagerungen von Kreisen, deren Radien nicht sehr
verschieden sind.} Technische Hochschule Stuttgart 1966. 42~Seiten.

\item
Heinrich W\"{o}lpert: \emph{Transformationstheorie der quadratischen Strahlkomplexe
der Charakteristik\/} $[(33)]$. Technische Hochschule Stuttgart 1967.

\item
Richard Koch: \emph{Geometrien mit einer Cayleyschen Fl\"{a}che dritten Grades als
absolutem Gebilde.} Universit\"{a}t Stuttgart (Technische Hochschule) 1968.
132~Seiten.

\item
Manfred Oehler: \emph{Axiomatisierung der Geometrie auf der Cayley'schen
Fl\"{a}che.} Universit\"{a}t Stuttgart (Technische Hochschule ) 1969. 71~Seiten.

\item
Lothar Profke: \emph{Kongruente Verlagerung projektiver Ebenen in Grenzlage.}
Universit\"{a}t Stuttgart 1969. 97~Seiten.

\item
Gunther R\"{o}sler: \emph{Zur Differentialgeometrie des Fl\"{a}chenelementes dritter
Ordnung.} Universit\"{a}t Stuttgart 1969. 44~Seiten.

\item
Siegfried Gr\"{u}ner: \emph{Zur Differentialgeometrie der isotropen M\"{o}biusebene.}
Universit\"{a}t Stuttgart 1970. 95~Seiten.

\item
Kurt Peter M\"{u}ller: \emph{Zur Geometrie der symplektischen Gruppe im reellen
dreidimensionalen projektiven Raum.} Universit\"{a}t Stuttgart 1970. 63~Seiten.

\item
Georg Kronhuber: \emph{Regelfl\"{a}chen im zentroaffinen Raum.} Technische
Hochschule Wien 1972. 52~Seiten.

\item
Gunter Wei{\ss}: \emph{\"{U}ber die Striktionslinie reeller analytischer windschiefer
Fl\"{a}chen.} Technische Hochschule Wien 1973.

\item
Rolf Riesinger: \emph{Verallgemeinerte zirkulare und verallgemeinerte zyklische
quadratische Komplexe.} Technische Hochschule Wien 1973. 114~Seiten.

\item
Hans-Peter Paukowitsch: \emph{Differentialgeometrie der Kurven bez\"{u}glich der
Scherungsgruppe in n-dimensionalen reellen affinen R\"{a}umen.} Technische
Hochschule Wien, 1975. 54~Seiten.

\item
Friedrich Anzb\"{o}ck: \emph{Eine durch das Kleinsche \"{U}bertragungsprinzip
vermittelte projektive Differentialgeometrie der windschiefen Fl\"{a}chen.}
Technische Hochschule Wien 1976, 87~Seiten.\\
Erschienen in: \emph{J.\ Reine Angew.\ Math.} \textbf{299}/\textbf{300} (1978),
92--112.

\item
Herbert Fritsche: \emph{Eine geometrische Kennzeichnung der linearen
Abbildungen der Geraden des projektiven, dreidimensionalen Raumes in eine
projektive Ebene.} Technische Universit\"{a}t Wien 1978. 56~Seiten.

\item
Hans Havlicek: \emph{Lineare Abbildungen aus Gra{\ss}mann-R\"{a}umen.}
Technische Universit\"{a}t Wien 1980. 100~Seiten.\\
Erschienen als: Zur Theorie linearer Abbildungen I, II. \emph{J.\ Geom.}
\textbf{16} (1981), 152--167, 168--180.

\item
Friedrich Manhart: \emph{Zur relativen Differentialgeometrie der Hyperfl\"{a}chen.}
Technische Universit\"{a}t Wien 1982. 68~Seiten.

\item
Ingrid Muhr: \emph{Fl\"{a}chen mit konstant geb\"{o}schten Schmieglinien.} Technische
Universit\"{a}t Wien 1983. 34~Seiten.

\item
Andreas Asperl: \emph{Zweifach Blutelsche Kegelschnittsfl\"{a}chen des projektiven
dreidimensionalen Raumes.} Technische Universit\"{a}t Wien 1990.
93~Seiten.\footnote{Andreas Asperl stellte seine von Heinrich Brauner betreute
Dissertation einen Monat nach dessen Ableben fertig. Erstgutachter war Gunter
Wei{\ss}.}

\end{publikationen}

\section*{Anhang 3: Dissertationen\\(Zweitgutachter H.~Brauner)}
\begin{publikationen}{9}

\item
Hubert Bitzel: \emph{Zur Konstruktion von \"{u}bertragungsg\"{u}nstigen, ebenen
Kurvengetrieben mit schwingendem oder umlaufendem Abtriebsglied.} Technische
Hochschule Stuttgart 1969. 80~Seiten. Erstgutachter: J.~Jehlicka.

\item
Gunther Petersch: \emph{Algebraische Raumkurven vierter Ordnung mit fester
Hauptnormalenneigung.} Technische Hochschule Wien 1974. 67~Seiten.
Erstgutachter: Walter Wunderlich.

\item
G\"{u}nter Eigenthaler: \emph{Zur Theorie der Polynome und Polynomfunktionen.}
Technische Hochschule Wien 1975. 93~Seiten. Erstgutachter: Winfried N\"{o}bauer.

\item
Maximilian Kreuzer: \emph{Eichunabh\"{a}ngige Schwelleneffekte bei
vereinheitlichten Theorien der starken, schwachen und elektromagnetischen
Wechselwirkungen.} Technische Universit\"{a}t Wien 1986, 109~Seiten. Erstgutachter:
Wolfgang Kummer.

\item
Hartwig Sorger: \emph{Eigenschattengrenzen konvexer K\"{o}rper und Verwandtes.}
Technische Universit\"{a}t Wien 1987, 68~Seiten. Erstgutachter: Peter M.~Gruber.

\item
Peter Grandits: \emph{C-Diskrepanz von Fl\"{a}chen im Raum.} Technische Universit\"{a}t
Wien 1990, 71~Seiten. Erstgutachter: Rudolf Taschner.

\end{publikationen}

\noindent Hans Havlicek, Forschungsgruppe Differentialgeometrie und
Geometrische Strukturen, Institut f\"{u}r Diskrete Mathematik und Geometrie,
Technische Universit\"{a}t Wien, Wiedner Hauptstra{\ss}e 8--10, A-1040 Wien,
Austria.
\newline
\texttt{havlicek@geometrie.tuwien.ac.at}

\end{document}